\documentclass[a4paper,12pt]{amsart}

\usepackage{amsmath,amssymb,amsthm,mathtools}
\usepackage{mathrsfs}
\usepackage[a4paper,margin=1in]{geometry}
\usepackage[
    colorlinks=true,
    linkcolor=blue,
    citecolor=magenta,
    urlcolor=blue
]{hyperref}
\usepackage{MnSymbol}
\mathtoolsset{showonlyrefs=true}

\newtheorem{theorem}{Theorem}[section]
\newtheorem{proposition}[theorem]{Proposition}
\newtheorem{lemma}[theorem]{Lemma}
\newtheorem{corollary}[theorem]{Corollary}

\theoremstyle{definition}
\newtheorem{definition}[theorem]{Definition}

\theoremstyle{remark}
\newtheorem{remark}[theorem]{Remark}

\numberwithin{equation}{section}

\newcommand{\bNzero}{\mathbb{N}_0}
\newcommand{\bZ}{\mathbb{Z}}
\newcommand{\bR}{\mathbb{R}}
\newcommand{\bC}{\mathbb{C}}
\newcommand{\bT}{\mathbb{T}}
\newcommand{\SD}{\mathscr{S}_D}
\newcommand{\Lip}{\operatorname{Lip}}
\newcommand{\BV}{\operatorname{BV}}
\newcommand{\Res}{\operatorname*{Res}}
\newcommand{\ind}{\mathbf{1}}
\newcommand{\dd}{\mathrm{d}}

\title[Periodic counting and carry synchronization]
{Periodic counting, carry synchronization, and primes
in nonstandard numeration systems}

\author{Nian Hong \textsc{Zhou}}

\address{
School of Mathematics and Statistics,
The Center for Applied Mathematics of Guangxi,
Guangxi Normal University,
Guilin, 541006, Guangxi, PR China
}

\email{nianhongzhou@outlook.com; nianhongzhou@gxnu.edu.cn}



\begin{document}

\begin{abstract}
Let $b\ge2$ be an integer and let $D\subset\mathbb{N}_0$ be a finite digit set containing zero, with maximum residue multiplicity $\mu=\max_{0\le r<b}\#\{d\in D:d\equiv r\pmod b\}$ and cardinality $q=\#D\ge2$. Under the condition $\mu<q$, we prove that the summatory representation function satisfies an asymptotic formula with a strictly positive, H\"older continuous, $1$-periodic fluctuating profile and an error of $O_D(x^{\log_b\mu})$. An exact periodic factorization of the Laplace generating function furnishes an explicit, computable Fourier expansion for the periodic profile via Gamma factors, and yields the meromorphic continuation of the associated Dirichlet series to the entire complex plane. When $D$ is a complete residue system modulo $b$, we establish an arithmetic synchronization criterion showing that the associated finite carry system admits a zero-synchronizing word if and only if $\gcd(D)=1$. This synchronization yields an exponentially accurate interval approximation and a general transfer principle for arithmetic weights. Consequently, we derive a prime number theorem for representable primes governed by the same periodic profile, alongside an exponentially convergent measurable scaling limit of the representable sets.
\end{abstract}

\maketitle

\section{Introduction}
\label{sec:introduction}

Fix an integer base $b\ge2$ and a finite digit set
$D\subset\bNzero$ containing zero.
The set of nonnegative integers representable in base $b$
with digits in $D$ is given by
\[
\SD
=
\left\{
\sum_{j\ge0}\delta_jb^j:
\delta_j\in D,\quad
\delta_j=0
\text{ for all sufficiently large }j
\right\}.
\]
When $D=\{0,1,\ldots,b-1\}$ is the canonical digit set,
every nonnegative integer admits a unique expansion and
$\SD=\bNzero$.
For a nonstandard digit set $D$, the arithmetic
structure of $\SD$ and its representation multiplicities
depend intricately on the magnitudes of the digits,
their distribution across residue classes modulo $b$,
and the carry propagation arising when converting from
standard base-$b$ expansions.

In this paper, we investigate representation counting,
the distribution of primes in $\SD$, and the explicit
harmonic analysis of the periodic fluctuations induced
by the multiplicative scaling $x\mapsto bx$.
Finite automata and nonstandard numeration systems
provide a natural framework for studying arithmetic
functions defined via digital expansions; see the
monograph of Allouche and Shallit \cite{AlloucheShallit2003}.
Periodic fluctuations in summatory functions governed by
divide-and-conquer recurrences, together with analytic
methods for computing their Fourier series, have been
systematically investigated by Dumas \cite{Dumas2014}
and by Heuberger and Krenn \cite{HeubergerKrenn2020}.

For $n\in\bNzero$, let $R_D(n)$ denote the number of
representations of $n$ in base $b$ with digits in $D$, namely
\[
R_D(n)
=
\#\left\{
(\delta_j)_{j\ge0}\in D^{\bNzero}:
n=\sum_{j\ge0}\delta_jb^j\in\SD
\right\},
\]
and define the summatory representation function
\[
A_D(x)=\sum_{0\le n\le x}R_D(n).
\]
Representations are regarded as eventually zero
sequences, so adjoining leading zeros does not
produce a new representation.

Put
\begin{equation}
\label{eq:parameters}
q=\#(D),
\qquad
\mu=
\max_{0\le r<b}
\#\{d\in D:d\equiv r\pmod b\},
\end{equation}
and define the exponents
\begin{equation}
\label{eq:exponents}
\alpha=\log_bq,
\qquad
\beta=\log_b\mu,
\qquad
\eta=\alpha-\beta.
\end{equation}
Since $q\le b\mu$, we always have $\eta\le1$.
Our standing assumption for the general counting theory
is $\mu<q$, which ensures that $\eta>0$.

Throughout the paper, the base $b$ is fixed.
Implied constants carrying a subscript $D$ may depend
on both $b$ and $D$.
For $0<\eta\le1$, we write $\Lip_\eta(\bR/\bZ)$ for the
space of $1$-periodic H\"older continuous functions
of exponent $\eta$, and $\BV(\bR/\bZ)$ for the space of
$1$-periodic functions of bounded variation on $[0,1]$.

\begin{theorem}[Periodic representation counting]
\label{thm:counting}
Suppose that $0\in D$, $q\ge2$, and $\mu<q$.
There exists a unique continuous $1$-periodic
function $\Psi_D$ such that
\begin{equation}
\label{eq:counting-asymptotic}
A_D(x)
=
x^\alpha\Psi_D(\log_bx)+O_D(x^\beta)
\qquad(x\ge1).
\end{equation}
The profile $\Psi_D$ is strictly positive and
satisfies
\[
\Psi_D\in
\Lip_\eta(\bR/\bZ)\cap\BV(\bR/\bZ).
\]
Moreover, there exists a constant $C_D>0$ such that the
one-sided error bound
\[
0\le
A_D(x)-x^\alpha\Psi_D(\log_bx)
\le C_Dx^\beta
\qquad(x\ge1)
\]
holds.
\end{theorem}

While $\Psi_D$ is constructed in Section~\ref{sec:counting}
as a monotone scaling limit of the step function $A_D$, we also
derive an explicit, rapidly computable Fourier expansion for
$\Psi_D$ that avoids evaluating $A_D(b^m x)$ at large scales $m$.
For $t>0$, let $P_D(z)=\sum_{d\in D}z^d$ and define the
infinite products
\[
F_D(t)
=
\sum_{n\ge0}R_D(n)e^{-nt}
=
\prod_{j\ge0}P_D\!\left(e^{-b^jt}\right),
\qquad
\mathcal{C}_D(t)
=
\prod_{j\ge1}
\left(q^{-1}P_D\!\left(e^{-t/b^j}\right)\right)^{-1}.
\]
Both products converge rapidly on compact subsets of $(0,\infty)$,
and the function 
$$v\mapsto b^{\alpha v}F_D(b^v)/\mathcal{C}_D(b^v)$$
is smooth and $1$-periodic on $\bR$.

\begin{theorem}[Explicit Fourier representation]
\label{thm:fourier-counting-profile}
Under the hypotheses of Theorem~\ref{thm:counting},
put $\tau_k=2\pi k/\log b$ and
\[
h_D(k)
=
\int_0^1
\frac{b^{\alpha v}F_D(b^v)}
{\mathcal{C}_D(b^v)}
e^{2\pi ikv}\,\dd v
\qquad(k\in\bZ).
\]
Then the Fourier coefficients of $\Psi_D$ are given by
\[
\widehat{\Psi}_D(k)
=
\int_0^1\Psi_D(u)e^{-2\pi iku}\,\dd u
=
\frac{h_D(k)}{\Gamma(1+\alpha+i\tau_k)},
\]
and the Fourier series of $\Psi_D$ converges absolutely
and uniformly on $\bR$:
\begin{equation}
\label{eq:counting-profile-gamma-fourier}
\Psi_D(u)
=
\sum_{k\in\bZ}
\frac{h_D(k)}
{\Gamma(1+\alpha+i\tau_k)}
e^{2\pi iku}.
\end{equation}
In particular, if $D$ is a complete residue system modulo $b$
containing zero, then
\begin{equation}
\label{eq:complete-residue-fourier}
\Psi_D(u)
=
\sum_{k\in\bZ}
\frac{h_D(k)}
{\Gamma(2+i\tau_k)}
e^{2\pi iku}.
\end{equation}
\end{theorem}

As a further application of this periodic factorization, we
show in Section~\ref{sec:dirichlet} that the Dirichlet
series $\zeta_D(s)=\sum_{n\ge1}R_D(n)n^{-s}$ admits a
meromorphic continuation to all of $\bC$, with simple
candidate poles at $s=\alpha-j+i\tau_k$ ($j\in\bNzero$, $k\in\bZ$)
whose residues are explicitly determined by $\widehat{\Psi}_D(k)$.
We note that Theorem~\ref{thm:counting} and Theorem \ref{thm:fourier-counting-profile}, and the Dirichlet series
continuation in Section~\ref{sec:dirichlet} specialize and
streamline results from the author's unpublished preprint
\cite{Zhou2023} on general real-base numeration systems; since
\cite{Zhou2023} is not intended for publication in that form, we
include self-contained and simplified proofs here for completeness.

\medskip

A particularly important regime arises when $D$ forms a
complete residue system modulo $b$:
\begin{equation}
\label{eq:complete-residue-assumption}
D=\{d_0,\ldots,d_{b-1}\},
\qquad
d_r\equiv r\pmod b,
\qquad
d_0=0.
\end{equation}
Successive determination of digits from least to most
significant shows that every integer admits at most one
finite $D$-representation, so that
\[
R_D(n)=\ind_{\SD}(n).
\]
In this case $q=b$ and $\mu=1$, and
Theorem~\ref{thm:counting} yields
\[
A_D(x)
=
\#(\SD\cap[0,x])
=
x\Psi_D(\log_bx)+O_D(1),
\]
where $\Psi_D$ is strictly positive and Lipschitz
continuous.

To study the distribution of primes in $\SD$ under
\eqref{eq:complete-residue-assumption}, we introduce a
finite carry automaton.
For an integer $t$, let $[t]_b\in\{0,\ldots,b-1\}$ denote its
least nonnegative residue modulo $b$.
For each standard input digit $a\in\{0,\ldots,b-1\}$, define
the carry transition map
\[
h_a(c)=b^{-1}(a+c-d_{[a+c]_b}).
\]
Setting
\[
d_*=\max D,
\qquad
M=\left\lceil\frac{d_*}{b-1}\right\rceil,
\]
one verifies that the finite state set
\[
K_D=\{-M,-M+1,\ldots,0\}
\]
is invariant under $h_a$ for every $a\in\{0,\ldots,b-1\}$.

Input words are read from the least significant
digit to the most significant digit.
For a word $w=a_0a_1\cdots a_{\ell-1}\in\{0,\ldots,b-1\}^\ell$, we write
\[
h_w=h_{a_{\ell-1}}\circ\cdots\circ h_{a_0}.
\]
We call $w$ \emph{zero-synchronizing} if
\[
h_w(K_D)=\{0\}.
\]

\begin{theorem}[Arithmetic criterion for synchronization]
\label{thm:synchronization}
Suppose that $D$ satisfies
\eqref{eq:complete-residue-assumption}.
The carry system on $K_D$ admits a
zero-synchronizing word if and only if
$\gcd(D)=1$.
\end{theorem}

The necessity of $\gcd(D)=1$ follows from a modular
congruence obstruction for differences of carry states.
To establish sufficiency, we study the induced Markov
chain on the finite set of carry differences.
If a nonempty closed subset of nonzero difference states
failed to reach zero, any stationary probability distribution
on that subset would generate a trigonometric polynomial
satisfying a weighted averaging identity on the unit circle.
Applying a minimum principle to the real part of this
trigonometric polynomial contradicts $\gcd(D)=1$.
Thus every pair of carry states is mergeable, and
iterative pairwise merging produces a zero-synchronizing word.

Carry synchronization leads directly to an interval
decomposition of $\SD$.
Specifically, a fixed zero-synchronizing word appears
in all but an exponentially small proportion of long
input blocks.
Whenever the high-order digital block of an integer contains
this synchronizing word, its terminal carry state---and hence
its membership in $\SD$---becomes completely independent of
all lower-order digits.
This allows uniform interval estimates for general
arithmetic weights to be transferred directly to $\SD$.

Prime numbers with restricted digital expansions
form a central theme in modern analytic number theory.
In breakthrough work, Maynard \cite{Maynard2019} proved that
for each $a_0\in\{0,\ldots,9\}$, the sparse set of positive
integers whose decimal expansions omit $a_0$ contains
$\asymp X^{\log 9/\log 10}/\log X$ primes up to $X$.
He subsequently established asymptotic formulas for
primes and prime values of polynomials with restricted
digits in sufficiently large bases \cite{Maynard2022},
combining the Hardy--Littlewood circle method with
Harman's sieve to decouple Diophantine and digital
Fourier conditions.
Our prime distribution theorem addresses a complementary
regime: for complete residue digit systems, $\SD$ has
positive upper and lower densities governed by the
periodic profile $\Psi_D$, while nonstandard digit sizes
induce unboundedly long carry chains when converting from
standard base-$b$ expansions.
Rather than bounding exponential sums over sparse
fractal sets of zero density, we harness carry
synchronization to approximate $\SD$ by unions of long
intervals, thereby transferring the classical prime number
theorem with error term.

Define the Chebyshev and prime counting functions restricted to $\SD$:
\[
\vartheta_D(x)
=
\sum_{\substack{p\le x\\p\in\SD}}\log p,
\qquad
\pi_D(x)
=
\#\{p\le x:p\text{ is prime and }p\in\SD\}.
\]

\begin{theorem}[Prime distribution]
\label{thm:primes}
Suppose that $D$ satisfies
\eqref{eq:complete-residue-assumption}
and that $\gcd(D)=1$.
There exists a constant $c_D>0$ such that
\begin{equation}
\label{eq:weighted-primes-intro}
\vartheta_D(x)
=
A_D(x)
+
O_D\!\left(xe^{-c_D\sqrt{\log x}}\right).
\end{equation}
In particular,
\begin{equation}
\label{eq:unweighted-primes-intro}
\pi_D(x)
=
\frac{x}{\log x}\Psi_D(\log_bx)
+
O_D\!\left(\frac{x}{(\log x)^2}\right).
\end{equation}
\end{theorem}

Thus, the exact same periodic profile $\Psi_D$ governs both
integer and prime counting in $\SD$.
Since $\Psi_D$ is strictly positive, we obtain the asymptotic
equivalence
\[
\pi_D(x)\sim\frac{A_D(x)}{\log x},
\]
and comparison with the classical prime number theorem yields
\[
\frac{\pi_D(x)}{\pi(x)}
=
\Psi_D(\log_bx)
+
O_D\!\left(\frac1{\log x}\right).
\]
Whenever $\Psi_D$ is nonconstant, the relative proportion of
primes lying in $\SD$ exhibits persistent logarithmic oscillations.

Carry synchronization also provides a geometric
interpretation of $\Psi_D$ in terms of a measurable
scaling limit.
For each integer $N\ge1$, define the rescaled subset of $[0,1)$:
\begin{equation}
\label{eq:scaled-sets-intro}
E_N=
\bigsqcup_{\substack{0\le n<b^N\\n\in\SD}}
\left[
\frac{n}{b^N},
\frac{n+1}{b^N}
\right).
\end{equation}

\begin{theorem}[Scaling limit]
\label{thm:scaling-limit}
Suppose that $D$ satisfies
\eqref{eq:complete-residue-assumption}
and that $\gcd(D)=1$.
There exist a measurable set $E_D\subset[0,1]$
and constants $C_D>0$ and $0<\rho_D<1$ such that
\[
\|\ind_{E_N}-\ind_{E_D}\|_{L^1([0,1])}
\le C_D\rho_D^N
\qquad(N\ge1).
\]
Moreover, for every $0<t\le1$,
\[
\int_0^t\ind_{E_D}(u)\,\dd u
=
t\Psi_D(\log_bt).
\]
\end{theorem}

\medskip

The remainder of the paper is organized as follows.
Section~\ref{sec:counting} establishes Theorem~\ref{thm:counting}.
Section~\ref{sec:dirichlet} derives the periodic factorization
of $F_D$, proves Theorem~\ref{thm:fourier-counting-profile},
and establishes the meromorphic continuation of $\zeta_D$.
Section~\ref{sec:carries} develops the finite carry automaton
and proves the synchronization criterion (Theorem~\ref{thm:synchronization}).
Section~\ref{sec:primes} establishes the cylinder interval
approximation and the arithmetic-weight transfer principle,
completing the proofs of Theorems~\ref{thm:primes} and \ref{thm:scaling-limit}.
Finally, Section~\ref{sec:examples} presents concrete examples
and discusses open problems.

\section{Representation counts and periodic asymptotics}
\label{sec:counting}

Throughout this section, $D\subset\bNzero$
is a finite digit set satisfying $0\in D$ and $q=|D|\ge2$,
with parameters defined as in \eqref{eq:parameters} and
\eqref{eq:exponents}.
For general algebraic and analytic approaches to
divide-and-conquer recurrences and $b$-regular
sequences, we refer to \cite{Dumas2014,HeubergerKrenn2020};
here we exploit the nonnegativity of the digit recurrence
to give a direct monotone scaling construction that
simultaneously yields one-sided error bounds, strict
positivity, and H\"older regularity.

\subsection{The digit recurrence}

For $0\le r<b$, put
\[
D_r=\{d\in D:d\equiv r\pmod b\}.
\]
We extend $R_D$ to all real arguments by setting
$R_D(y)=0$ if $y\notin\bNzero$.
Since the digits in $D$ are nonnegative, the only
representation of $0$ is the all-zero sequence,
so $R_D(0)=1$.

\begin{lemma}[Digit recurrence]
\label{lem:digit-recurrence}
For every $n\in\bNzero$,
\[
R_D(n)
=
\sum_{d\in D}
R_D\!\left(\frac{n-d}{b}\right)
=
\sum_{d\in D_{[n]_b}}
R_D\!\left(\frac{n-d}{b}\right).
\]
\end{lemma}

\begin{proof}
Every representation of $n$ can be decomposed uniquely as
\[
n=\delta_0+
b\sum_{j\ge0}\delta_{j+1}b^j.
\]
Once the least significant digit $\delta_0=d\in D$ is fixed,
the remaining digits $(\delta_{j+1})_{j\ge0}$ form a representation
of $(n-d)/b$, which is an integer if and only if $d\equiv n\pmod b$.
\end{proof}

\begin{proposition}[Pointwise bound]
\label{prop:pointwise-bound}
For every integer $n\ge1$,
\[
R_D(n)
\le
\mu^{1+\lfloor\log_bn\rfloor}
\le
\mu n^\beta.
\]
\end{proposition}

\begin{proof}
Since every nonzero digit in $D$ is at least $1$, a nonzero
digit in position $j$ contributes at least $b^j$.
Thus every representation of $n\ge1$ uses at most
$m=1+\lfloor\log_bn\rfloor$ positions.
Padding all representations with leading zeros to length $m$,
we see from Lemma~\ref{lem:digit-recurrence} that at each of
the $m$ steps there are at most $\mu=\max_r|D_r|$ admissible
choices for the current digit.
Hence $R_D(n)\le\mu^m\le\mu^{1+\log_bn}=\mu n^\beta$.
\end{proof}

In particular, when $D$ is a complete residue system modulo $b$,
we have $\mu=1$ and $R_D(n)=\ind_{\SD}(n)$.

Setting $A_D(x)=0$ for $x<0$ and summing
Lemma~\ref{lem:digit-recurrence} over $0\le n\le x$ yields
the exact functional recurrence
\begin{equation}
\label{eq:summatory-recurrence}
A_D(x)
=
\sum_{d\in D}
A_D\!\left(\frac{x-d}{b}\right).
\end{equation}

\subsection{A monotone scaling limit}

For $x>0$ and $m\in\bNzero$, define the rescaled counting
sequence
\[
G_m(x)=q^{-m}A_D(b^mx).
\]
By \eqref{eq:summatory-recurrence} and the monotonicity of $A_D$,
\[
A_D(b^{m+1}x)
=
\sum_{d\in D}
A_D\!\left(b^mx-\frac{d}{b}\right)
\le qA_D(b^mx).
\]
Consequently, $G_{m+1}(x)\le G_m(x)$ for all $m\ge0$, so the
nonnegative pointwise limit
\[
\mathcal{F}_D(x)
=
\lim_{m\to\infty}q^{-m}A_D(b^mx)
\]
exists for all $x>0$.

\begin{lemma}[Successive scales]
\label{lem:successive-scales}
There exists a constant $C_D>0$ such that
\[
0\le qA_D(y)-A_D(by)\le C_Dy^\beta
\qquad(y\ge1).
\]
\end{lemma}

\begin{proof}
By \eqref{eq:summatory-recurrence},
\[
qA_D(y)-A_D(by)
=
\sum_{d\in D}
\left(
A_D(y)-A_D\!\left(y-\frac{d}{b}\right)
\right).
\]
Each summand on the right-hand side sums $R_D(n)$ over integers
$n$ in an interval of bounded length $d/b\le d_*/b$.
For $n\ge1$, Proposition~\ref{prop:pointwise-bound} gives
$R_D(n)\le\mu y^\beta$, and the possible term $R_D(0)=1$
satisfies the same bound since $y\ge1$.
\end{proof}

Assume now that $\mu<q$.
Dividing the estimate of Lemma~\ref{lem:successive-scales}
at $y=b^m x$ ($x\ge1$) by $q^{m+1}$ gives
\[
0\le G_m(x)-G_{m+1}(x)
\le
\frac{C_D}{q}x^\beta
\left(\frac{\mu}{q}\right)^m.
\]
Summing this telescoping bound over $m\ge0$ (with $C_D$
redefined by a constant factor $1/(q-\mu)$) yields
\[
0\le G_m(x)-\mathcal{F}_D(x)
\le
C_Dx^\beta
\left(\frac{\mu}{q}\right)^m.
\]
In particular, taking $m=0$ gives
\begin{equation}
\label{eq:counting-limit-error}
0\le A_D(x)-\mathcal{F}_D(x)
\le C_Dx^\beta
\qquad(x\ge1).
\end{equation}

By construction, $\mathcal{F}_D$ satisfies the scaling relation
$\mathcal{F}_D(bx)=q\mathcal{F}_D(x)$.
Define
\[
\Psi_D(u)
=
b^{-\alpha u}\mathcal{F}_D(b^u),
\qquad u\in\bR.
\]
Since $b^\alpha=q$, the function $\Psi_D$ is $1$-periodic on $\bR$,
and $\mathcal{F}_D(x)=x^\alpha\Psi_D(\log_bx)$.
Thus \eqref{eq:counting-limit-error} establishes both the asymptotic
formula \eqref{eq:counting-asymptotic} and the one-sided error bound
in Theorem~\ref{thm:counting}.
Moreover, uniformly for $u\in[0,1]$ and $m\ge0$,
\[
0\le
b^{-\alpha u}q^{-m}A_D(b^{m+u})
-\Psi_D(u)
\le
C_D\left(\frac{\mu}{q}\right)^m.
\]

\subsection{Short intervals and regularity}

For $\ell\in\bNzero$, let
\[
R_{D,\ell}(a)
=
\#\left\{
(\delta_0,\ldots,\delta_{\ell-1})\in D^\ell:
a=\sum_{j=0}^{\ell-1}\delta_jb^j
\right\},
\]
with the convention $R_{D,0}(0)=1$ and $R_{D,0}(a)=0$ for $a\ne0$.
Splitting any representation at position $\ell$ gives the
convolution identity
\begin{equation}
\label{eq:block-decomposition}
R_D(n)
=
\sum_{\substack{a,k\ge0\\a+b^\ell k=n}}
R_{D,\ell}(a)R_D(k).
\end{equation}

\begin{lemma}[Relative short-interval estimate]
\label{lem:short-interval}
Fix $C>0$.
There exists a constant $C_{D,C}>0$ such that, for all
integers $N\ge m\ge0$ and every interval
$I\subset[0,Cb^N]$ of length $|I|\le b^{N-m}$,
one has
\[
\sum_{n\in I\cap\bNzero}R_D(n)
\le C_{D,C}q^{N-m}\mu^m.
\]
\end{lemma}

\begin{proof}
Set $\ell=N-m$ and $B_D=d_*/(b-1)$, where $d_*=\max D$.
Any lower block of length $\ell$ has value
$0\le a\le B_D(b^\ell-1)\le B_Db^\ell$.
Thus, if $n=a+b^\ell k\in I$, then $k$ must lie in the interval
$b^{-\ell}I-[0,B_D]$, which has length at most $1+B_D$ and is
contained in $[0,Cb^m]$.
Hence at most $2+\lfloor B_D\rfloor$ integers $k$ can contribute,
and for each such $k$, Proposition~\ref{prop:pointwise-bound}
and $R_D(0)=1$ give $R_D(k)\le C_{D,C}'\mu^m$.
Since $\sum_{a\ge0}R_{D,\ell}(a)=q^\ell=q^{N-m}$, summing
\eqref{eq:block-decomposition} over $n\in I\cap\bNzero$
yields the desired bound.
\end{proof}

To establish the H\"older continuity of $\Psi_D$, let $1\le x<y\le b$
and suppose first that $h=y-x\le1$.
Choose the integer $m\ge0$ such that $b^{-m-1}<h\le b^{-m}$.
For any $N\ge m$, the interval $(b^Nx,b^Ny]\subset[0,b^{N+1}]$
has length $b^N h\le b^{N-m}$.
Applying Lemma~\ref{lem:short-interval} with $C=b$ gives
\[
q^{-N}
\bigl(A_D(b^Ny)-A_D(b^Nx)\bigr)
\le
C_{D,b}\left(\frac{\mu}{q}\right)^m.
\]
Letting $N\to\infty$ and using $\mu/q=b^{-\eta}$ and $b^{-m}<bh$,
we obtain
\[
0\le
\mathcal{F}_D(y)-\mathcal{F}_D(x)
\le C_D'|y-x|^\eta.
\]
For $h=y-x>1$, the same estimate holds trivially by the
boundedness of $\mathcal{F}_D$ on $[1,b]$.

Since the maps $u\mapsto b^u$ and $u\mapsto b^{-\alpha u}$ are
smooth on $[0,1]$, it follows that
\[
|\Psi_D(u)-\Psi_D(v)|
\le C_D''|u-v|^\eta
\qquad(0\le u,v\le1).
\]
Because $\Psi_D(0)=\Psi_D(1)$, applying the triangle
inequality across the period endpoints shows that the
$1$-periodic extension preserves this H\"older bound, so
$\Psi_D\in\Lip_\eta(\bR/\bZ)$.

\begin{proposition}[Positivity and bounded variation]
\label{prop:positive-bv}
The function $\Psi_D$ is strictly positive,
bounded away from zero, and belongs to
$\BV(\bR/\bZ)$.
\end{proposition}

\begin{proof}
Put $B_D=d_*/(b-1)>0$.
Every digit string of length $m$ represents an integer
strictly smaller than $B_Db^m$, and there are $q^m$ such strings
counted with representation multiplicity.
Hence $A_D(B_Db^m)\ge q^m$ for all $m\ge0$, which implies
$\mathcal{F}_D(B_D)\ge1$.

For any $x>0$, choose an integer $j\ge0$ such that $b^jx\ge B_D$.
By monotonicity and the scaling relation,
\[
q^j\mathcal{F}_D(x)
=
\mathcal{F}_D(b^jx)
\ge\mathcal{F}_D(B_D)\ge1,
\]
so $\mathcal{F}_D(x)>0$.
By continuity and $1$-periodicity, $\Psi_D$ attains a strictly
positive minimum on $\bR$.

Finally, the function $u\mapsto\mathcal{F}_D(b^u)$ is bounded and
nondecreasing on $[0,1]$, and hence of bounded variation.
Since multiplication by the smooth function $b^{-\alpha u}$
preserves bounded variation and $\Psi_D(0)=\Psi_D(1)$, we
conclude that $\Psi_D\in\BV(\bR/\bZ)$.
\end{proof}

\begin{proof}[Completion of the proof of Theorem~\ref{thm:counting}]
It remains only to verify the uniqueness of $\Psi_D$.
Suppose that $\widetilde{\Psi}_D$ is another continuous $1$-periodic
function satisfying \eqref{eq:counting-asymptotic}.
For any fixed $t>0$, substituting $x=b^mt$ and dividing by
$(b^mt)^\alpha$ gives
\[
\Psi_D(\log_bt)
-
\widetilde{\Psi}_D(\log_bt)
=
O_D\!\left(
t^{\beta-\alpha}b^{-m(\alpha-\beta)}
\right).
\]
Since $\alpha>\beta$, letting $m\to\infty$ yields
$\Psi_D(\log_bt)=\widetilde{\Psi}_D(\log_bt)$.
\end{proof}

\begin{corollary}[Complete residue systems]
\label{cor:complete-residue-counting}
If $D$ is a complete residue system modulo $b$ containing zero,
there exists a strictly positive Lipschitz continuous $1$-periodic
function $\Psi_D$ such that
\[
\#(\SD\cap[0,x])
=
x\Psi_D(\log_bx)+O_D(1)
\qquad(x\ge1).
\]
\end{corollary}

\begin{proof}
Here $q=b$, $\mu=1$, $\eta=1$, and $R_D(n)=\ind_{\SD}(n)$.
\end{proof}

\begin{remark}
If $\gcd(D)=1$, the condition $\mu<q$ holds automatically.
Indeed, if $\mu=q$, all elements of $D$ lie in a single residue
class modulo $b$; since $0\in D$, every digit in $D$ would be
divisible by $b$, contradicting $\gcd(D)=1$.
Note also that Corollary~\ref{cor:complete-residue-counting}
requires no coprimality assumption; the condition $\gcd(D)=1$
enters only in the synchronization and prime distribution theorems.
\end{remark}

\section{Explicit Fourier expansion and Dirichlet series}
\label{sec:dirichlet}

In this section, we return to the general setting of a finite digit
set $D\subset\bNzero$ satisfying
\[
0\in D,\qquad q=|D|\ge2,\qquad\mu<q,
\]
without any completeness or synchronization hypothesis.

\subsection[Periodic factorization]{An exact periodic factorization}

Define the digit polynomial and ordinary generating function
\[
P_D(z)=\sum_{d\in D}z^d,
\qquad
Z_D(z)=\sum_{n\ge0}R_D(n)z^n
=
\prod_{j\ge0}P_D(z^{b^j})
\qquad(|z|<1).
\]
For $t>0$, the Laplace generating function
\[
F_D(t)
=
Z_D(e^{-t})
=
\sum_{n\ge0}R_D(n)e^{-nt}
\]
satisfies the functional equation $F_D(t)=q\,p_D(t)F_D(bt)$, where
$p_D(t)=P_D(e^{-t})/q$.

Since $p_D(0)=1$, we may choose a radius $\delta_D>0$ such that
$|p_D(t)-1|<1/2$ for all $t\in\bC$ with $|t|<\delta_D$, and define
the principal analytic logarithm $L_D(t)=\log p_D(t)$ with $L_D(0)=0$.
Because $b\ge2$ and $L_D(t)=O_D(t)$ on $|t|<\delta_D$, the series
$\sum_{j\ge1}L_D(t/b^j)$ converges uniformly on compact subsets of
$|t|<\delta_D$, so
\[
\mathcal{C}_D(t)
=
\exp\left(
-\sum_{j\ge1}L_D(t/b^j)
\right)
=
\prod_{j\ge1}\bigl(p_D(t/b^j)\bigr)^{-1}
\]
defines a nonvanishing analytic function on $|t|<\delta_D$ with
$\mathcal{C}_D(0)=1$, and extends as a smooth, strictly positive
function to all real $t>0$ satisfying
\[
\mathcal{C}_D(bt)
=
\frac{\mathcal{C}_D(t)}{p_D(t)}.
\]
We now define
\[
H_D(u)
=
\frac{b^{\alpha u}F_D(b^u)}
{\mathcal{C}_D(b^u)},
\qquad u\in\bR.
\]

\begin{proposition}[Periodic factorization]
\label{prop:periodic-factorization}
The function $H_D$ is strictly positive, smooth, and $1$-periodic
on $\bR$, and for every $t>0$,
\[
F_D(t)
=
t^{-\alpha}
\mathcal{C}_D(t)H_D(\log_bt).
\]
Moreover, $\mathcal{C}_D(t)$ admits a Taylor expansion near $t=0$
of the form
\[
\mathcal{C}_D(t)
=
\sum_{j\ge0}c_D(j)t^j,
\qquad c_D(0)=1.
\]
\end{proposition}

\begin{proof}
Combining $F_D(b^u)=q\,p_D(b^u)F_D(b^{u+1})$,
$\mathcal{C}_D(b^{u+1})=\mathcal{C}_D(b^u)/p_D(b^u)$, and
$b^\alpha=q$ immediately gives $H_D(u+1)=H_D(u)$.
Since $F_D$ and $\mathcal{C}_D$ are both smooth and strictly
positive on $(0,\infty)$, so is $H_D$, and the Taylor expansion
of $\mathcal{C}_D$ around $t=0$ follows from its analyticity on
$|t|<\delta_D$.
\end{proof}

\subsection{Fourier coefficients and the Gamma relation}

Recall that
\[
\widehat{\Psi}_D(k)
=
\int_0^1\Psi_D(u)e^{-2\pi iku}\,\dd u,
\qquad
\tau_k=\frac{2\pi k}{\log b},
\]
and define the Fourier coefficients of $H_D$ with the opposite
sign convention:
\[
h_D(k)
=
\int_0^1H_D(u)e^{2\pi iku}\,\dd u,
\qquad
H_D(u)
=
\sum_{k\in\bZ}h_D(k)e^{-2\pi iku}.
\]
Since $H_D\in C^\infty(\bR/\bZ)$, repeated integration by parts gives
\begin{equation}
\label{eq:fourier-rapid-decay}
h_D(k)
=
O_{D,A}\bigl((1+|k|)^{-A}\bigr)
\end{equation}
for every $A\ge0$.

\begin{proposition}[Gamma relation]
\label{prop:gamma-relation}
For every $k\in\bZ$,
\[
h_D(k)
=
\Gamma(1+\alpha+i\tau_k)
\widehat{\Psi}_D(k).
\]
\end{proposition}

\begin{proof}
By Laplace--Stieltjes integration by parts,
\[
F_D(t)
=
t\int_0^\infty A_D(x)e^{-tx}\,\dd x
\qquad(t>0).
\]
Substituting Theorem~\ref{thm:counting} on $[1,\infty)$ and
changing variables to $v=tx$, we find as $t\to0^+$,
\[
t^\alpha F_D(t)
=
\int_0^\infty
v^\alpha
\Psi_D(\log_bv-\log_bt)e^{-v}\,\dd v
+
O_D(t^\eta),
\]
where the error $O_D(x^\beta)$ on $[1,\infty)$ contributes
$O_D(t^{\alpha-\beta})=O_D(t^\eta)$ and the integral over $0<x<1$
contributes $O_D(t^{\alpha+1})$.

Fixing $u\in\bR$ and setting $t=b^{u-m}$ with $m\to\infty$, we deduce
from the $1$-periodicity of $H_D$ and $\Psi_D$ together with
$\lim_{t\to0^+}\mathcal{C}_D(t)=1$ that
\begin{equation}
\label{eq:profile-integral-transform}
H_D(u)
=
\int_0^\infty
v^\alpha
\Psi_D(\log_bv-u)e^{-v}\,\dd v.
\end{equation}
Multiplying both sides of \eqref{eq:profile-integral-transform} by
$e^{2\pi iku}$, integrating over $u\in[0,1]$, and applying Fubini's
theorem (since $\Psi_D$ is bounded), we obtain by $1$-periodicity
\[
\int_0^1
\Psi_D(\log_bv-u)e^{2\pi iku}\,\dd u
=
e^{2\pi ik\log_bv}\widehat{\Psi}_D(k)
=
v^{i\tau_k}\widehat{\Psi}_D(k),
\]
and hence
\[
h_D(k)
=
\widehat{\Psi}_D(k)
\int_0^\infty
v^{\alpha+i\tau_k}e^{-v}\,\dd v
=
\Gamma(1+\alpha+i\tau_k)\widehat{\Psi}_D(k).
\qedhere
\]
\end{proof}

In particular, for $k=0$,
\[
h_D(0)
=
\Gamma(1+\alpha)
\int_0^1\Psi_D(u)\,\dd u
>0.
\]

\subsection{Proof of the explicit Fourier representation}

\begin{proof}[Proof of Theorem~\ref{thm:fourier-counting-profile}]
By Theorem~\ref{thm:counting}, $\Psi_D$ belongs to
$\Lip_\eta(\bR/\bZ)\cap\BV(\bR/\bZ)$ with $\eta=\alpha-\beta>0$.
By Zygmund's theorem on functions of bounded variation satisfying a
H\"older condition \cite[Chapter~VI, Theorem~(3.6), p.~241]{Zygmund2002},
the Fourier series of $\Psi_D$ converges absolutely and uniformly
on $\bR$.
Combining this with Proposition~\ref{prop:gamma-relation} proves
\eqref{eq:counting-profile-gamma-fourier}.

In the special case of a complete residue system modulo $b$
(\eqref{eq:complete-residue-fourier}), we have $\alpha=1$ and
$\Psi_D\in\Lip_1(\bR/\bZ)$, so absolute convergence also follows
directly from Rademacher's theorem ($\Psi_D'\in L^\infty([0,1])\subset L^2([0,1])$),
Parseval's identity $\sum_{k\in\bZ}k^2|\widehat{\Psi}_D(k)|^2<\infty$,
and the Cauchy--Schwarz inequality.
\end{proof}

\subsection{Mellin continuation and residues}

The correspondence between logarithmically periodic fluctuations
and vertical arithmetic progressions of complex poles is a hallmark
of Mellin asymptotics; see the survey of Flajolet, Gourdon, and
Dumas \cite{FlajoletGourdonDumas1995}, as well as
\cite{Dumas2014,HeubergerKrenn2020} for the explicit link between
residues of Dirichlet generating functions and Fourier coefficients
of periodic profiles.
As a consequence of Propositions~\ref{prop:periodic-factorization}
and \ref{prop:gamma-relation}, we establish the meromorphic
continuation of the Dirichlet series
\[
\zeta_D(s)
=
\sum_{n\ge1}\frac{R_D(n)}{n^s}.
\]

By Theorem~\ref{thm:counting}, $A_D(x)\asymp_D x^\alpha$, so partial
summation shows that $\zeta_D(s)$ converges absolutely for
$\operatorname{Re}s>\alpha$ and diverges at $s=\alpha$.
For $\operatorname{Re}s>\alpha$, Mellin inversion gives
\[
\Gamma(s)\zeta_D(s)
=
\int_0^\infty
t^{s-1}\bigl(F_D(t)-1\bigr)\,\dd t.
\]
Fix $\varepsilon\in(0,\delta_D)$ so that $\mathcal{C}_D$ is analytic
on a neighborhood of $[0,\varepsilon]$.
The tail integral $\int_\varepsilon^\infty t^{s-1}(F_D(t)-1)\,\dd t$
defines an entire function of $s\in\bC$.

For any integer $J\ge0$, expanding $\mathcal{C}_D(t)$ to order $J$
via Proposition~\ref{prop:periodic-factorization} gives
\[
F_D(t)
=
\sum_{j=0}^J
c_D(j)t^{j-\alpha}H_D(\log_bt)
+
O_{D,J}(t^{J+1-\alpha})
\qquad(0<t\le\varepsilon),
\]
where the Mellin transform of the remainder over $(0,\varepsilon)$
is holomorphic on the half-plane $\operatorname{Re}s>\alpha-J-1$.
Since $\sum_{k\in\bZ}|h_D(k)|<\infty$, Fubini's theorem justifies
termwise integration of the Fourier expansion
$H_D(\log_bt)=\sum_{k\in\bZ}h_D(k)t^{-i\tau_k}$ for
$\operatorname{Re}s>\alpha-j$, yielding
\[
\int_0^\varepsilon
t^{s+j-\alpha-1}H_D(\log_bt)\,\dd t
=
\sum_{k\in\bZ}
h_D(k)
\frac{
\varepsilon^{s+j-\alpha-i\tau_k}
}{
s+j-\alpha-i\tau_k
}.
\]
By the rapid decay \eqref{eq:fourier-rapid-decay} of $h_D(k)$, the
series on the right-hand side defines a meromorphic function on all
of $\bC$.
Thus, for $\operatorname{Re}s>\alpha-J-1$,
\begin{align}
\Gamma(s)\zeta_D(s)
={}&
\sum_{j=0}^Jc_D(j)
\sum_{k\in\bZ}
h_D(k)
\frac{
\varepsilon^{s+j-\alpha-i\tau_k}
}{
s+j-\alpha-i\tau_k
}-
\frac{\varepsilon^s}{s}
+
\mathcal{R}_{D,J}(s),
\label{eq:continued-mellin-expression}
\end{align}
where $\mathcal{R}_{D,J}(s)$ is holomorphic on
$\operatorname{Re}s>\alpha-J-1$.

\begin{theorem}[Pole locations and residues]
\label{thm:full-meromorphic-continuation}
The Dirichlet series $\zeta_D(s)$ extends meromorphically to $\bC$.
Its candidate poles are all simple and lie on the vertical lines
\[
s_{j,k}=\alpha-j+i\tau_k,
\qquad
j\in\bNzero,\quad k\in\bZ.
\]
At each candidate point $s_{j,k}$, the residue is given by
\begin{equation}
\label{eq:general-residue-formula}
\Res_{s=s_{j,k}}\zeta_D(s)
=
\frac{c_D(j)h_D(k)}
{\Gamma(\alpha-j+i\tau_k)}
=
c_D(j)
\frac{\Gamma(1+\alpha+i\tau_k)}
{\Gamma(\alpha-j+i\tau_k)}
\widehat{\Psi}_D(k),
\end{equation}
where $1/\Gamma(s)$ is understood as an entire function vanishing at
nonpositive integers.
In particular, on the abscissa of convergence $\operatorname{Re}s=\alpha$,
\[
\Res_{s=\alpha+i\tau_k}\zeta_D(s)
=
(\alpha+i\tau_k)\widehat{\Psi}_D(k).
\]
\end{theorem}

\begin{proof}
Letting $J\to\infty$ in \eqref{eq:continued-mellin-expression}
continues $\Gamma(s)\zeta_D(s)$ meromorphically to $\bC$ with
simple candidate poles at $s=0$ and at the distinct points $s_{j,k}$.
Dividing by $\Gamma(s)$ introduces no new singularities; the simple
zero of $1/\Gamma(s)$ at $s=0$ cancels the pole of $-\varepsilon^s/s$,
and its simple zeros at $s\in-\bNzero$ cancel any simple poles of
the double sum occurring at nonpositive integers.
Applying Proposition~\ref{prop:gamma-relation} gives
\eqref{eq:general-residue-formula}, and setting $j=0$ with
$c_D(0)=1$ yields the simplified formula on $\operatorname{Re}s=\alpha$.
\end{proof}

In particular, $s=\alpha$ is always a genuine simple pole with
strictly positive residue
\[
\Res_{s=\alpha}\zeta_D(s)
=
\alpha
\int_0^1\Psi_D(u)\,\dd u
>0.
\]

\begin{remark}[Constant profiles and rightmost poles]
\label{rem:pole-criterion}
The periodic profile $\Psi_D$ is constant if and only if all
nonreal candidate poles $s_{0,k}=\alpha+i\tau_k$ ($k\ne0$) on
$\operatorname{Re}s=\alpha$ are removable, since $\alpha+i\tau_k\ne0$
and $\widehat{\Psi}_D(k)=0$ for all $k\ne0$ if and only if $\Psi_D$
is constant.
\end{remark}

\section{Finite carry systems and synchronization}
\label{sec:carries}

Throughout this section, we assume that
\begin{equation}
\label{eq:carry-digit-assumption}
D=\{d_0,\ldots,d_{b-1}\}\subset\bNzero,
\qquad
d_r\equiv r\pmod b,
\qquad
d_0=0.
\end{equation}
We write $\mathcal{A}=\{0,1,\ldots,b-1\}$ for the standard base-$b$
alphabet and extend the indexing of $D$ periodically to all $r\in\bZ$
by setting $d_r=d_{[r]_b}$.

In the terminology of Lagarias and Wang \cite{LagariasWang1996},
a digit set $D$ forming a complete residue system modulo $b$ is a
standard digit set for the expanding map $x\mapsto bx$.
In \cite{LagariasWang1996}, Fourier-analytic criteria involving the
mask polynomial $b^{-1}\sum_{d\in D}e^{2\pi idy}$ were developed
to study the Lebesgue measure of self-affine tiles $T(b,D)$.
The weighted averaging identity established in
Lemma~\ref{lem:weighted-averaging} below uses the squared modulus
of this trigonometric polynomial to analyze reachability in the
finite difference system of carries.

\subsection{Carry maps and terminal states}

For $a\in\mathcal{A}$ and $c\in\bZ$, recall the carry transition map
\[
h_a(c)=b^{-1}(a+c-d_{a+c}),
\]
which satisfies the fundamental identity
\[
a+c=d_{a+c}+bh_a(c).
\]
Put
\[
M=\left\lceil\frac{d_*}{b-1}\right\rceil,
\qquad
K_D=\{-M,-M+1,\ldots,0\}.
\]

\begin{lemma}[Invariant state set]
\label{lem:invariant-state-set}
For every $a\in\mathcal{A}$, one has $h_a(K_D)\subseteq K_D$.
\end{lemma}

\begin{proof}
For $c\in K_D$, since $c\le0$ and $d_{a+c}\ge0$, we have
$h_a(c)\le a/b\le(b-1)/b<1$.
Because $h_a(c)\in\bZ$, it follows that $h_a(c)\le0$.
On the other hand, since $a\ge0$, $c\ge-M$, and $d_{a+c}\le d_*\le(b-1)M$,
we have
\[
h_a(c)\ge\frac{-M-d_*}{b}\ge-M.
\qedhere
\]
\end{proof}

For a finite input word $w=a_0\cdots a_{m-1}\in\mathcal{A}^m$,
read from least to most significant position, we set
\[
h_w=h_{a_{m-1}}\circ\cdots\circ h_{a_0},
\]
with the empty word acting as the identity map.
Thus concatenation of words satisfies $h_{uv}=h_v\circ h_u$.

\begin{proposition}[Terminal-state criterion]
\label{prop:terminal-state}
Let $0\le n<b^m$ have the standard base-$b$ expansion
\[
n=\sum_{j=0}^{m-1}a_jb^j.
\]
Starting from the initial carry $c_0=0$, define recursively
$c_{j+1}=h_{a_j}(c_j)$ for $0\le j<m$.
Then
\[
n\in\SD
\quad\Longleftrightarrow\quad
c_m=0.
\]
\end{proposition}

\begin{proof}
Set $\delta_j=d_{a_j+c_j}\in D$.
By the carry identity, $a_j+c_j=\delta_j+bc_{j+1}$.
Multiplying by $b^j$ and summing over $0\le j\le m-1$ gives
the telescoping relation
\begin{equation}
\label{eq:telescoping-carries}
n=\sum_{j=0}^{m-1}\delta_jb^j+b^mc_m.
\end{equation}
Thus, if $c_m=0$, then $n=\sum_{j=0}^{m-1}\delta_jb^j\in\SD$.

Conversely, suppose that $n\in\SD$.
Since every nonzero digit in $D$ is at least $1$ and $n<b^m$,
all digits of $n$ in positions $m$ and above must vanish.
Successive reduction modulo $b,b^2,\ldots,b^m$ forces the
digits of $n$ in positions $0\le j<m$ to coincide with the digits
$\delta_j$ generated by the carry recurrence, and
\eqref{eq:telescoping-carries} then forces $c_m=0$.
\end{proof}

\begin{definition}
A finite word $w\in\mathcal{A}^*$ is \emph{zero-synchronizing} if
$h_w(K_D)=\{0\}$.
Two states $c,c'\in K_D$ are \emph{mergeable} if there exists a
finite word $u\in\mathcal{A}^*$ such that $h_u(c)=h_u(c')$.
\end{definition}

\subsection{The arithmetic obstruction and pairwise merging}

\begin{proposition}[Necessary arithmetic condition]
\label{prop:gcd-necessary}
If a zero-synchronizing word exists, then $\gcd(D)=1$.
\end{proposition}

\begin{proof}
Let $g=\gcd(D)$.
Since every $d\in D$ is divisible by $g$, the definition of $h_a$
implies that for any $a\in\mathcal{A}$ and $c,c'\in K_D$,
\[
b\bigl(h_a(c)-h_a(c')\bigr)
=
(c-c')-\bigl(d_{a+c}-d_{a+c'}\bigr)
\equiv c-c'\pmod g.
\]
Iterating over a word $w$ of length $m$ yields
\[
b^m\bigl(h_w(c)-h_w(c')\bigr)
\equiv c-c'\pmod g.
\]
If $w$ is zero-synchronizing, then $h_w(c)=h_w(c')=0$, so $g$
must divide $c-c'$ for all $c,c'\in K_D$.
Since $|D|=b\ge2$, we have $d_*\ge b-1$ and $M\ge1$, so
$0,-1\in K_D$.
Taking $c=0$ and $c'=-1$ gives $g=1$.
\end{proof}

\begin{remark}
If $g=\gcd(D)>1$, then $\SD\subseteq g\bNzero$.
Consequently, any prime in $\SD$ must equal $g$, so $\SD$
contains at most one prime.
\end{remark}

\begin{lemma}[Accessibility of zero]
\label{lem:zero-accessible}
For every $c\in K_D$, there exists a finite word $v\in\mathcal{A}^*$
such that $h_v(c)=0$.
\end{lemma}

\begin{proof}
For $c=0$, the empty word suffices.
For $c<0$, write the standard base-$b$ expansion of $-c>0$ as
\[
-c=\sum_{j=0}^{m-1}a_jb^j.
\]
Reading the word $v=a_0\cdots a_{m-1}$ from $c_0=c$, we observe
that $a_0+c_0\equiv0\pmod b$.
Since $d_0=0$, we obtain
\[
c_1
=
\frac{a_0+c_0}{b}
=
-\sum_{j=1}^{m-1}a_jb^{j-1}.
\]
Repeating this step inductively yields $c_m=0$.
\end{proof}

\begin{proposition}[Pairwise merging criterion]
\label{prop:pairwise-criterion}
A zero-synchronizing word exists if and only if every pair of
states in $K_D$ is mergeable.
\end{proposition}

\begin{proof}
The forward implication is immediate.
Conversely, suppose that every pair of states in $K_D$ is mergeable.
Whenever the current image of $K_D$ contains at least two distinct
states, we may append a word that merges two of them, strictly
decreasing the cardinality of the image.
After at most $|K_D|-1$ such steps, we obtain a word $u$ such that
$h_u(K_D)=\{c_*\}$ for some $c_*\in K_D$.
Appending a word that sends $c_*$ to $0$ via
Lemma~\ref{lem:zero-accessible} produces a zero-synchronizing word.
\end{proof}

\subsection{The difference system}

All differences between carry states lie in the symmetric interval
\[
\mathcal{J}_D=K_D-K_D=\{-M,\ldots,M\}.
\]
For $r,k\in\bZ$, define
\[
g_r(k)=\frac{k+d_r-d_{r+k}}{b}.
\]
Since $d_{r+k}\equiv r+k\equiv k+d_r\pmod b$, each $g_r$ maps
$\bZ$ into $\bZ$ and satisfies $g_{r+b}=g_r$.
Direct calculation shows that for any $a\in\mathcal{A}$ and $c,c'\in K_D$,
\[
h_a(c)-h_a(c')
=
g_{a+c'}(c-c').
\]
Moreover, for $k\in\mathcal{J}_D$,
\[
|g_r(k)|
\le\frac{M+d_*}{b}\le M,
\]
so $\mathcal{J}_D$ is invariant under every $g_r$.

We define a Markov transition matrix $A=(A_{k,\ell})_{k,\ell\in\mathcal{J}_D}$ by
\[
A_{k,\ell}
=
\frac1b
\#\{r\in\{0,\ldots,b-1\}:g_r(k)=\ell\}.
\]
Clearly, $A$ is row-stochastic and $0\in\mathcal{J}_D$ is an
absorbing state ($A_{0,0}=1$).
Since $r\mapsto d_r$ bijects $\{0,\ldots,b-1\}$ onto $D$, we may
equivalently write
\[
A_{k,\ell}
=
\frac1b
\#\{(u,v)\in D^2:k+u-v=b\ell\}.
\]

\subsection{A weighted averaging identity}

Suppose toward a contradiction that some difference state cannot
reach $0$, and consider the nonempty set of such states:
\[
U=
\{k\in\mathcal{J}_D:
(A^n)_{k,0}=0\text{ for all }n\ge0\}.
\]
Since $A^0$ is the identity matrix, $0\notin U$.
If $k\in U$ and $A_{k,\ell}>0$, then necessarily $\ell\in U$
(otherwise a path from $\ell$ to $0$ would yield a path from $k$ to $0$).
Thus the restricted matrix $A_U=(A_{k,\ell})_{k,\ell\in U}$ is itself
row-stochastic, and hence admits a stationary probability vector
$(p_k)_{k\in U}$ satisfying
\[
p_k\ge0,\qquad
\sum_{k\in U}p_k=1,\qquad
\sum_{k\in U}p_kA_{k,\ell}=p_\ell
\quad(\ell\in U).
\]

On the torus $\bT=\bR/\bZ$, define the trigonometric polynomial
\[
H(x)=\sum_{k\in U}p_ke^{2\pi ikx},
\qquad
f(x)=\operatorname{Re}H(x),
\]
and the nonnegative weight function
\[
\omega(y)
=
\frac1{b^2}
\left|\sum_{d\in D}e^{2\pi idy}\right|^2.
\]

\begin{lemma}[Weighted averaging identity]
\label{lem:weighted-averaging}
For every $x\in\bR$,
\begin{equation}
\label{eq:weights-partition}
\sum_{r=0}^{b-1}
\omega\left(\frac{x+r}{b}\right)=1,
\end{equation}
and
\begin{equation}
\label{eq:weighted-harmonic}
H(x)
=
\sum_{r=0}^{b-1}
\omega\left(\frac{x+r}{b}\right)
H\left(\frac{x+r}{b}\right).
\end{equation}
\end{lemma}

\begin{proof}
Expanding the squared modulus gives
\[
\omega(y)
=
\frac1{b^2}
\sum_{u,v\in D}e^{2\pi i(u-v)y}.
\]
By orthogonality of $b$-th roots of unity,
summing $\omega((x+r)/b)$ over $0\le r<b$ retains only those
pairs $(u,v)\in D^2$ with $u\equiv v\pmod b$.
Since $D$ is a complete residue system modulo $b$, these are
precisely the diagonal pairs $u=v$, yielding \eqref{eq:weights-partition}.

Similarly,
\begin{align*}
\sum_{r=0}^{b-1}
\omega\left(\frac{x+r}{b}\right)
H\left(\frac{x+r}{b}\right)&=
\frac1b
\sum_{k\in U}p_k
\sum_{\substack{u,v\in D\\b\mid k+u-v}}
e^{2\pi i(k+u-v)x/b}\\
&=
\sum_{\ell\in\mathcal{J}_D}
\left(
\sum_{k\in U}p_kA_{k,\ell}
\right)e^{2\pi i\ell x}\\
&=
\sum_{\ell\in U}p_\ell e^{2\pi i\ell x}
=H(x),
\end{align*}
where the last line follows from the closure of $U$ and the
stationarity of $(p_k)_{k\in U}$.
\end{proof}

\subsection{The minimum principle}

\begin{lemma}[Zero-difference reachability]
\label{lem:difference-reaches-zero}
If $\gcd(D)=1$, every state in $\mathcal{J}_D$ can reach $0$
through a finite sequence of maps $g_r$.
\end{lemma}

\begin{proof}
Assume toward a contradiction that $U$ is nonempty.
Since $0\notin U$, we have
\[
f(0)=\sum_{k\in U}p_k=1,
\qquad
\int_0^1f(x)\,\dd x=0.
\]
Hence $f$ is nonconstant, and its global minimum
$m_f=\min_{x\in\bT}f(x)$ is strictly negative.

Consider the minimum set
\[
E=\{x\in[0,1):f(x)=m_f\}.
\]
Clearly $E$ is nonempty, $0\notin E$, and $E$ is finite because
$f-m_f$ is a nonzero trigonometric polynomial.

Taking real parts in \eqref{eq:weighted-harmonic} shows that for
every $x\in E$,
\[
m_f
=
\sum_{r=0}^{b-1}
\omega\left(\frac{x+r}{b}\right)
f\left(\frac{x+r}{b}\right).
\]
Since the weights $\omega((x+r)/b)$ are nonnegative and sum to $1$
by \eqref{eq:weights-partition}, every branch $(x+r)/b$ carrying
strictly positive weight must itself attain the minimum $m_f$.
Thus, setting
\[
P(x)=
\left\{
\frac{x+r}{b}:0\le r<b,\quad
\omega\left(\frac{x+r}{b}\right)>0
\right\},
\]
we have $\varnothing\ne P(x)\subseteq E$ for every $x\in E$.
Moreover, the sets $P(x)$ for distinct $x\in E$ are pairwise
disjoint, since $y\in P(x)$ uniquely determines $x=\{by\}$.
Consequently,
\[
\#E
\ge
\sum_{x\in E}\#P(x)
\ge
\#E,
\]
which forces $\#P(x)=1$ for every $x\in E$.
By \eqref{eq:weights-partition}, the single element $y\in P(x)$
must satisfy $\omega(y)=1$.

By the condition for equality in the triangle inequality and $0\in D$,
\[
\omega(y)=1
\quad\Longleftrightarrow\quad
e^{2\pi idy}=1
\text{ for all }d\in D
\quad\Longleftrightarrow\quad
dy\in\bZ
\text{ for all }d\in D.
\]
Since $\gcd(D)=1$, B\'ezout's identity implies $y\in\bZ$, so
$y=0$ is the unique point in $[0,1)$ with $\omega(y)=1$.
Thus $P(x)=\{0\}$ for any $x\in E$.
Since $P(x)\subseteq E$, this yields $0\in E$, contradicting
$0\notin E$.
Therefore $U=\varnothing$.
\end{proof}

\subsection{Construction of a synchronizing word}

\begin{lemma}[Pairwise coalescence]
\label{lem:pairwise-coalescence}
If $\gcd(D)=1$, every pair of states in $K_D$ is mergeable.
\end{lemma}

\begin{proof}
Given $c,c'\in K_D$ with $c\ne c'$, Lemma~\ref{lem:difference-reaches-zero}
provides a finite sequence of residues $r_0,\ldots,r_{n-1}$ such that
\[
k_0=c-c',
\qquad
k_{i+1}=g_{r_i}(k_i),
\qquad
k_n=0.
\]
Set $c_0=c$ and $c'_0=c'$.
At each step $0\le i<n$, choose the unique digit $a_i\in\mathcal{A}$
satisfying $a_i+c'_i\equiv r_i\pmod b$, and set
$c_{i+1}=h_{a_i}(c_i)$ and $c'_{i+1}=h_{a_i}(c'_i)$.
Then
\[
c_{i+1}-c'_{i+1}
=
g_{a_i+c'_i}(c_i-c'_i)
=
g_{r_i}(k_i)
=
k_{i+1}.
\]
Hence $c_n=c'_n$, so the word $a_0\cdots a_{n-1}$ merges $c$ and $c'$.
\end{proof}

\begin{proof}[Proof of Theorem~\ref{thm:synchronization}]
Necessity is given by Proposition~\ref{prop:gcd-necessary}.
Sufficiency follows immediately by combining
Lemma~\ref{lem:pairwise-coalescence} with
Proposition~\ref{prop:pairwise-criterion}.
Since $h_0(0)=0$, the resulting zero-synchronizing word may always
be chosen nonempty.
\end{proof}

\begin{remark}[Constructive algorithm]
The proof of Theorem~\ref{thm:synchronization} furnishes an explicit
finite algorithm to construct a zero-synchronizing word: one finds
shortest paths to $0$ in the finite directed graph on $\mathcal{J}_D$,
translates them into pair-merging words via Lemma~\ref{lem:pairwise-coalescence}
to iteratively collapse $K_D$ to a singleton, and finally appends
the word from Lemma~\ref{lem:zero-accessible}.
\end{remark}

\section{Interval approximation and arithmetic weights}
\label{sec:primes}

Throughout this section, we assume that $D$ satisfies
\eqref{eq:carry-digit-assumption} and $\gcd(D)=1$, and we fix a
nonempty zero-synchronizing word
\[
w=w_0\cdots w_{\ell-1}\in\mathcal{A}^\ell,
\qquad \ell\ge1.
\]
The principle that all but an exponentially small proportion of
words of length $r$ contain a fixed synchronizing word---so that
the state of a finite automaton is determined solely by the
letters following the last occurrence of that word---was employed
by Deshouillers, Drmota, and M\"ullner
\cite{DeshouillersDrmotaMuellner2015} to prove that automatic
sequences generated by synchronizing automata satisfy the Sarnak
conjecture.
In our setting, because digits are processed from least to most
significant, an occurrence of $w$ inside the high-order digital
block renders the terminal carry state independent of all
lower-order digits, thereby decomposing $\SD$ into a union of
discrete intervals up to an exponentially small exceptional set.

\subsection{Determined and exceptional cylinders}

For integers $N\ge r\ge1$ and $0\le u<b^r$, define the discrete
cylinder
\[
\mathcal{C}_{N,r}(u)
=
[ub^{N-r},(u+1)b^{N-r})\cap\bNzero.
\]
As $u$ ranges over $0\le u<b^r$, the cylinders $\mathcal{C}_{N,r}(u)$
partition $[0,b^N)\cap\bNzero$ into $b^r$ blocks of length $b^{N-r}$.
Expanding
\[
u=\sum_{j=0}^{r-1}a_jb^j,
\]
we associate to $u$ its length-$r$ input word $a_0\cdots a_{r-1}\in\mathcal{A}^r$,
read from low to high positions.

\begin{lemma}[Constancy on synchronized cylinders]
\label{lem:cylinder-constancy}
If the length-$r$ input word of $u$ contains $w$ as a consecutive
subword, then for every $n\in\mathcal{C}_{N,r}(u)$,
\[
n\in\SD
\quad\Longleftrightarrow\quad
u\in\SD.
\]
\end{lemma}

\begin{proof}
Write $n=y+ub^{N-r}$ with $0\le y<b^{N-r}$.
The length-$N$ input word of $n$ is the length-$(N-r)$ input word
of $y$ followed by the length-$r$ input word of $u$.
Factor the latter as $v_1wv_2$.
Comparing the carry trajectories of $n$ and $u$ starting from $0$,
we see that their carry states just before reading $w$ both belong
to $K_D$.
Since $h_w(K_D)=\{0\}$, both trajectories are reset to $0$
immediately after reading $w$, and hence reach the exact same
terminal state after reading the common suffix $v_2$.
The claim now follows from Proposition~\ref{prop:terminal-state}.
\end{proof}

We call a cylinder $\mathcal{C}_{N,r}(u)$ \emph{determined} if
the length-$r$ word of $u$ contains $w$, and \emph{exceptional}
otherwise.

\begin{lemma}[Words avoiding a fixed word]
\label{lem:avoiding-words}
The proportion of words in $\mathcal{A}^r$ avoiding $w$ as a
consecutive subword is at most $(1-b^{-\ell})^{\lfloor r/\ell\rfloor}$.
Consequently, there exist constants $C_D>0$ and $0<\rho_D<1$ such
that at most $C_Db^r\rho_D^r$ words of length $r$ avoid $w$.
\end{lemma}

\begin{proof}
Partitioning a length-$r$ word into $k=\lfloor r/\ell\rfloor$
disjoint blocks of length $\ell$ plus a remainder of length
$r-k\ell$, we see that any word avoiding $w$ cannot match $w$ on
any of the $k$ blocks.
Hence there are at most $(b^\ell-1)^kb^{r-k\ell}=b^r(1-b^{-\ell})^k$
such words, and the bound follows with $\rho_D=(1-b^{-\ell})^{1/\ell}$.
\end{proof}

\begin{proposition}[Interval approximation]
\label{prop:interval-approximation}
For $N\ge r\ge1$, every determined cylinder $\mathcal{C}_{N,r}(u)$
is either entirely contained in $\SD$ or entirely disjoint from $\SD$,
while the union of all exceptional cylinders contains at most
$C_Db^N\rho_D^r$ integers.
\end{proposition}

\begin{proof}
This follows immediately from Lemmas~\ref{lem:cylinder-constancy}
and \ref{lem:avoiding-words}, since each of the at most
$C_Db^r\rho_D^r$ exceptional cylinders contains $b^{N-r}$ integers.
\end{proof}

Letting $\mathcal{G}_{N,r}$ denote the union of all determined
representable cylinders and $\mathcal{B}_{N,r}$ the union of all
exceptional cylinders, we have
\[
\mathcal{G}_{N,r}
\subseteq\SD\cap[0,b^N)
\subseteq
\mathcal{G}_{N,r}\cup\mathcal{B}_{N,r},
\qquad
\frac{
\#\bigl(
(\SD\cap[0,b^N))
\triangle\mathcal{G}_{N,r}
\bigr)
}{b^N}
\le C_D\rho_D^r.
\]

\subsection[Arithmetic-weight transfer]{Transfer to general arithmetic weights}
\label{subsec:arithmetic-weight-transfer}

The cylinder decomposition cleanly separates the digital
combinatorics of $\SD$ from the distribution of general
arithmetic weights on intervals.

\begin{proposition}[Arithmetic-weight transfer]
\label{prop:arithmetic-weight-transfer}
Let $a:\bNzero\to\bC$ be an arithmetic function.
For $X=b^N$, suppose that there exist quantities
$m(X)\in\bC$, $E(X)\ge0$, and $W(X)\ge0$ such that,
uniformly for all half-open intervals $I\subseteq[0,X)$,
\begin{equation}
\label{eq:weight-interval-estimate}
\left|
\sum_{n\in I\cap\bNzero}a(n)
-
m(X)\#(I\cap\bNzero)
\right|
\le E(X),
\end{equation}
and
\begin{equation}
\label{eq:weight-pointwise-bound}
|a(n)|+|m(X)|\le W(X)
\qquad(0\le n<X).
\end{equation}
Then, for every $1\le r\le N$, uniformly for all real $0\le y<X$,
\begin{align*}
\sum_{\substack{0\le n\le y\\n\in\SD}}a(n)
-
m(X)A_D(y)=
O_D\!\left(
b^rE(X)+X\rho_D^rW(X)
\right).
\end{align*}
\end{proposition}

\begin{proof}
Intersect the cylinders $\mathcal{C}_{N,r}(u)$ ($0\le u<b^r$)
with $[0,\lfloor y\rfloor+1)$.
Each nonempty intersection consists of the integer points of a
half-open interval contained in $[0,X)$.
Since membership in $\SD$ is constant on each determined cylinder,
only the determined representable pieces and the exceptional pieces
contribute to the difference
\[
\sum_{\substack{0\le n\le y\\n\in\SD}}
\bigl(a(n)-m(X)\bigr).
\]
There are at most $b^r$ determined representable pieces, and by
\eqref{eq:weight-interval-estimate} their combined contribution is
bounded by $b^rE(X)$.
By Proposition~\ref{prop:interval-approximation}, the exceptional
pieces contain at most $C_DX\rho_D^r$ integers in total, so
\eqref{eq:weight-pointwise-bound} bounds their contribution by
$C_DX\rho_D^rW(X)$.
Finally, since $A_D(y)=\#(\SD\cap[0,y])$, the stated bound follows.
\end{proof}

The resolution parameter $r$ balances the interval sum error
against the exceptional set bound.
For instance, whenever $E(X)\ll Xe^{-a\sqrt{\log X}}$ and
$W(X)\ll(\log X)^B$ for fixed $a>0$ and $B\ge0$, choosing
$r=\lfloor\kappa\sqrt{\log X}\rfloor$ with $\kappa>0$ sufficiently
small yields
\[
\sum_{\substack{0\le n\le y\\n\in\SD}}a(n)
=
m(X)A_D(y)
+
O_D\!\left(
Xe^{-c\sqrt{\log X}}
\right)
\]
uniformly for $0\le y<X$, for some constant $c>0$.

\subsection{The prime number theorem}

Let $\vartheta(y)=\sum_{p\le y}\log p$.
By the classical prime number theorem with the de la Vall\'ee Poussin
error term (see Davenport \cite[Chapter~18]{Davenport2000}), there
exists an absolute constant $a>0$ such that
\begin{equation}
\label{eq:classical-pnt}
\vartheta(y)
=
y+O\!\left(
ye^{-a\sqrt{\log y}}
\right),
\qquad y\ge2.
\end{equation}

\begin{lemma}[Uniform interval estimate]
\label{lem:uniform-endpoints}
There exists an absolute constant $a_1>0$ such that, for all
sufficiently large $X$, uniformly for $0\le u\le v\le X$,
\[
\sum_{u\le p<v}\log p
=
v-u+
O\!\left(
Xe^{-a_1\sqrt{\log X}}
\right).
\]
\end{lemma}

\begin{proof}
For $X^{1/2}\le y\le X$, \eqref{eq:classical-pnt} gives
$\vartheta(y)-y=O(Xe^{-(a/\sqrt2)\sqrt{\log X}})$, while for
$0\le y<X^{1/2}$ the trivial bound $\vartheta(y)\le y\log(2+y)$
is absorbed into the same error term upon choosing $0<a_1<a/\sqrt2$.
Subtracting the estimates at $y=v$ and $y=u$ and absorbing the
endpoint contribution $O(\log X)$ yields the lemma.
\end{proof}

\begin{proof}[Proof of Theorem~\ref{thm:primes}]
Given $x\ge2$, choose the integer $N\ge1$ such that
$b^{N-1}\le x<b^N$, and set $X=b^N\asymp_b x$.
Applying Proposition~\ref{prop:arithmetic-weight-transfer} to the
von Mangoldt prime weight $a(n)=\ind_{\text{$n$ prime}}\log n$,
we see from Lemma~\ref{lem:uniform-endpoints} that
\eqref{eq:weight-interval-estimate} and \eqref{eq:weight-pointwise-bound}
hold with
\[
m(X)=1,\qquad
E(X)\ll Xe^{-a_1\sqrt{\log X}},
\qquad
W(X)\ll\log X.
\]
Consequently, for any $1\le r\le N$,
\[
\vartheta_D(x)-A_D(x)
=
O_D\!\left(
b^rXe^{-a_1\sqrt{\log X}}
+
X\rho_D^r\log X
\right).
\]
Choosing $r=\lfloor\kappa\sqrt{\log X}\rfloor$ with
$0<\kappa<a_1/(2\log b)$ balances the two error terms and proves
\eqref{eq:weighted-primes-intro}.

Combining \eqref{eq:weighted-primes-intro} with
Corollary~\ref{cor:complete-residue-counting} gives
\[
\vartheta_D(x)
=
x\Psi_D(\log_bx)
+
O_D\!\left(
xe^{-c_D\sqrt{\log x}}
\right).
\]
By partial summation,
\[
\pi_D(x)
=
\frac{\vartheta_D(x)}{\log x}
+
\int_2^x
\frac{\vartheta_D(t)}
{t(\log t)^2}\,\dd t.
\]
Since $0\le\vartheta_D(t)\le\vartheta(t)\ll t$, the integral on the
right-hand side is $O(x/(\log x)^2)$, which establishes
\eqref{eq:unweighted-primes-intro}.
\end{proof}

\subsection{The measurable scaling limit}

Recall the rescaled sets $E_N\subset[0,1)$ defined in
\eqref{eq:scaled-sets-intro}, whose indicator functions satisfy
$\ind_{E_N}(t)=\ind_{\SD}(\lfloor b^Nt\rfloor)$ for $t\in[0,1)$.

\begin{proof}[Proof of Theorem~\ref{thm:scaling-limit}]
Fix integers $m\ge N\ge1$.
For $t\in[0,1)$ outside the countable set of $b$-adic rationals,
write its base-$b$ expansion as $t=0.a_1a_2\cdots$.
Then $n_N=\lfloor b^Nt\rfloor$ and $n_m=\lfloor b^mt\rfloor$
share the exact same $N$ most significant digits $a_1,\ldots,a_N$.

When testing $n_N$, the carry automaton reads the word
$a_Na_{N-1}\cdots a_1$ starting from $0$; when testing $n_m$, it
first processes the $m-N$ lower-order digits $a_m\cdots a_{N+1}$
to reach some carry state in $K_D$, and then reads the common word
$a_Na_{N-1}\cdots a_1$.
If $a_Na_{N-1}\cdots a_1$ contains the synchronizing word $w$, both
trajectories reach the same terminal state, so
$\ind_{E_m}(t)=\ind_{E_N}(t)$.

Since the prefix $(a_1,\ldots,a_N)$ is uniformly distributed on
$\mathcal{A}^N$ with respect to Lebesgue measure on $[0,1)$,
Lemma~\ref{lem:avoiding-words} implies that the exceptional set of
$t\in[0,1)$ where $\ind_{E_m}(t)\ne\ind_{E_N}(t)$ has measure at most
$C_D\rho_D^N$.
Thus
\begin{equation}
\label{eq:scaled-set-cauchy}
\|\ind_{E_m}-\ind_{E_N}\|_{L^1([0,1])}
\le C_D\rho_D^N
\qquad(m\ge N\ge1).
\end{equation}
Hence $(\ind_{E_N})_{N\ge1}$ is a Cauchy sequence in $L^1([0,1])$
with exponential rate $\rho_D^N$, and its $L^1$-limit is the
indicator function $\ind_{E_D}$ of a measurable set $E_D\subset[0,1]$.
Letting $m\to\infty$ in \eqref{eq:scaled-set-cauchy} gives the
stated $L^1$ error bound.

Finally, for any fixed $0<t\le1$,
\[
\int_0^t\ind_{E_N}(u)\,\dd u
=
b^{-N}A_D(b^Nt)+O(b^{-N})
=
t\Psi_D(N+\log_bt)+O_D(b^{-N})
\]
by Corollary~\ref{cor:complete-residue-counting}.
Using the $1$-periodicity of $\Psi_D$ and letting $N\to\infty$
completes the proof.
\end{proof}

\section{Examples and further questions}
\label{sec:examples}

\subsection{A nonconstant periodic profile}

Consider the base $b=3$ and the complete residue digit set
$D=\{0,1,5\}$, which satisfies $\gcd(D)=1$.
Here $d_*=5$, $M=\lceil 5/2\rceil=3$, $K_D=\{-3,-2,-1,0\}$,
and the carry transition table is given by
\[
\begin{array}{c|rrr}
c & h_0(c) & h_1(c) & h_2(c)\\
\hline
-3 & -1 & -1 & -2\\
-2 & -1 & -2 & 0\\
-1 & -2 & 0 & 0\\
0 & 0 & 0 & -1
\end{array}
\]
Tracing images of $K_D$, we find
\[
h_0(K_D)=\{-2,-1,0\},
\qquad
h_2(h_0(K_D))=\{-1,0\},
\qquad
h_1(\{-1,0\})=\{0\},
\]
so the low-to-high input word $w=021$ is zero-synchronizing.

By Corollary~\ref{cor:complete-residue-counting} and
Theorem~\ref{thm:primes}, there exists a strictly positive
Lipschitz continuous $1$-periodic profile $\Psi_D$ such that
\[
A_D(x)
=
x\Psi_D(\log_3x)+O_D(1),
\qquad
\pi_D(x)
=
\frac{x}{\log x}\Psi_D(\log_3x)
+
O_D\!\left(
\frac{x}{(\log x)^2}
\right),
\]
and by Theorem~\ref{thm:fourier-counting-profile}, $\Psi_D$ admits
the absolutely and uniformly convergent Fourier series
\[
\Psi_D(u)
=
\sum_{k\in\bZ}
\frac{h_D(k)}
{\Gamma\!\left(
2+\frac{2\pi ik}{\log3}
\right)}
e^{2\pi iku}.
\]

\subsection{A missing-digit example}

Consider the ternary Cantor digit set $b=3$, $D=\{0,1\}$, for which
\[
q=2,\qquad
\mu=1,\qquad
\alpha=\log_32,\qquad
\beta=0.
\]
Here representations are unique and Theorem~\ref{thm:counting} gives
\[
A_D(x)
=
x^{\log_32}\Psi_D(\log_3x)+O_D(1),
\]
with $\Psi_D\in\Lip_{\log_32}(\bR/\bZ)\cap\BV(\bR/\bZ)$ strictly
positive.
By Theorem~\ref{thm:fourier-counting-profile}, its Fourier expansion
\[
\Psi_D(u)
=
\sum_{k\in\bZ}
\frac{h_D(k)}
{\Gamma\!\left(
1+\log_32+\frac{2\pi ik}{\log3}
\right)}
e^{2\pi iku}
\]
converges absolutely and uniformly on $\bR$, and $\zeta_D(s)$
extends meromorphically to $\bC$ with simple candidate poles at
$\log_32-j+2\pi ik/\log3$ ($j\in\bNzero$, $k\in\bZ$).
While Theorem~\ref{thm:primes} does not apply here because $D$ is
not a complete residue system modulo $3$ (for prime counting in
sparse missing-digit sets, see \cite{Maynard2019,Maynard2022}),
this example illustrates that our counting and Fourier-analytic
results hold for a substantially broader class of digit sets than
the synchronization and prime theorems.

\subsection{Quantitative synchronization and arithmetic weights}

Theorem~\ref{thm:synchronization} settles the qualitative
existence problem: a zero-synchronizing word exists on $K_D$ if and
only if $\gcd(D)=1$.
A natural quantitative problem is to establish upper bounds on the
minimal length $\ell(D)$ of a zero-synchronizing word in terms of
$b$ and $d_*=\max D$ (in the spirit of \v{C}ern\'y's conjecture for
synchronizing automata).

Similarly, while Lemma~\ref{lem:avoiding-words} bounds the number
of $w$-avoiding words by partitioning into disjoint blocks of
length $\ell$, sharper exponential decay rates $\rho_D$ can be
obtained either from the spectral radius of the subshift transition
matrix avoiding $w$, or directly from the subset automaton on
$\mathcal{P}(K_D)\setminus\{\{c\}:c\in K_D\}$.

Finally, beyond prime indicator weights, the transfer principle of
Proposition~\ref{prop:arithmetic-weight-transfer} applies immediately
to any arithmetic function $a(n)$ admitting uniform cancellation on
intervals, such as the M\"obius function $\mu(n)$ or $n\mapsto\mu(n)e(n\theta)$,
yielding M\"obius orthogonality for the indicator function
$\ind_{\SD}$ in complete analogy with the Sarnak conjecture for
synchronizing automatic sequences \cite{DeshouillersDrmotaMuellner2015}.

\subsection*{Acknowledgements}
This work was partially supported by the Guangxi Natural Science Foundation (No.~2026GXNSFBA00640038) and the National Natural Science Foundation of China (No.~12661065). The authors used OpenAI's ChatGPT for exploratory discussions concerning possible proof strategies and language editing. All mathematical arguments, verification, and responsibility for the final manuscript strictly remain with the authors.

\end{document}